\documentclass{amsart}
\usepackage{amsmath,amssymb,amscd,graphicx,psfrag,epsf,amsthm
}
\usepackage[all]{xy}

\newcommand{\ou}{{\mathcal {O}}}

\newcommand{\cL}{{\mathcal {L}}}

\newcommand{\nef}{\overline{NM}^1}
\newcommand{\eff}{\overline{NE}^1}
\newcommand{\NEe}{\overline{NE}_1(X)_{(K_X+\Delta+ A)\geq 0}}
\newcommand{\NE}{\overline{NE}_1(X)_{(K_X+\Delta)\geq 0}}

\newcommand{\NEb}{\overline{NE}_1(X)_{K_X\geq 0}}
\newcommand{\nm}{\overline{NM}_1(X)}
\numberwithin{equation}{section}

\newtheorem*{thm*}{Theorem}

\newtheorem{thm}{Theorem}[section]

\newtheorem{lemma}[thm]{Lemma}

\newtheorem{cor}[thm]{Corollary}

\theoremstyle{definition}
\newtheorem{defn}[thm]{Definition}

\newtheorem{say}[thm]{}

\newtheorem{exmp}[thm]{Example}

\newtheorem{rem}[thm]{Remark}

\newtheorem{defn-thm}[thm]{Definition-Theorem}  

\theoremstyle{remark}

\newcommand{\z}[0]{{\mathbb Z}}

\renewcommand{\r}[0]{{\mathbb R}}

\newcommand{\p}[0]{{\mathbb P}}

\newcommand{\q}[0]{{\mathbb Q}}
\newcommand{\map}[0]{\dasharrow}

\newcommand{\pic}[0]{\operatorname{Pic}}

\newcommand{\im}[0]{\operatorname{im}}

\newcommand{\inter}[0]{\operatorname{Int}}

\newcommand{\nec}[0]{\overline{NE}_1}

\begin{document}

\title[The cone of pseudo-effective divisors]{The cone of pseudo-effective divisors of log varieties after Batyrev}
\author{Carolina Araujo}
\address{Carolina Araujo\newline\indent 
  IMPA -  Estrada Dona Castorina 110,
  Rio de Janeiro, Brazil 22460-320}
\email{caraujo@impa.br}

\begin{abstract}
In these notes we investigate the cone of nef curves of projective varieties, 
which is the dual cone to the cone of pseudo-effective divisors.
We prove a structure theorem for  the cone of nef curves of projective
$\q$-factorial klt pairs  of arbitrary dimension from the point of view of the Minimal Model Program. 
This is a generalization of Batyrev's structure theorem for the cone of nef curves of projective terminal threefolds.
\end{abstract}

\maketitle


\section{Introduction}

Cones of curves and divisors play  important roles in the classification of complex projective varieties,
specially in the Minimal Model Program. 
The goal of the Minimal Model Program is to produce, starting with a  $\q$-factorial klt  pair $(X,\Delta)$,
a birational map $\varphi:X\map X'$ such that either
$(X',\varphi_*\Delta)$ is a minimal model, or $X'$ admits a Mori Fiber Space structure.
We refer to  section~\ref{section_mmp} for definitions and results concerning the Minimal Model Program.
Unless otherwise stated, we always allow $\Delta$ to have real coefficients.
When $\Delta=\emptyset$ and $X$ is a terminal threefold, the program was established in \cite{mori_mmp}.
Recently, a special instance of the Minimal Model Program, the  Minimal Model Program with scaling,
has been established for certain types of pairs in arbitrary dimension in \cite{bchm}.
As a consequence, for any $\q$-factorial klt  pair $(X,\Delta)$  such that $K_X+\Delta$ is not pseudo-effective,
i.e., such that  $K_X+\Delta$ is not a limit of effective $\q$-divisors, 
one can run  a Minimal Model Program  to produce a 
birational map $\varphi:X\map X'$ and  a Mori Fiber Space structure on $X'$.

Given  a normal projective variety $X$, we denote by $N_1(X)$ the finite dimensional $\r$-vector space
of $1$-cycles on $X$ with real coefficients modulo numerical equivalence.
We denote by $\overline{NE}_1(X)$ the closed convex cone in $N_1(X)$ generated by  classes
of irreducible curves on $X$. 
One has the following structure theorem.
For a proof, see \cite[Theorem 3.7]{kollar_mori} when $\Delta$ has rational coefficients, 
and \cite[Corollary 3.10]{shokurov_3fold_log_models} in the general case. 
Here and throughout the paper, given a Cartier divisor $D$ on $X$,
we denote by $\overline{NE}_1(X)_{D\geq 0}$ the subcone
$\big\{z\in \overline{NE}_1(X) \ \big|\ D\cdot z\geq 0\big\}$. 

\begin{thm*} [Cone Theorem]
Let $(X,\Delta)$ be a klt  pair. There is a set $\Gamma\subset \overline{NE}_1(X)$
of classes of rational curves $C\subset X$ with $0<-(K_X+\Delta)\cdot C \leq 2\dim(X)$
such that 
\begin{enumerate}
\item for any ample $\r$-divisor $A$ on $X$,
		there are finitely many classes $[C_1],\dots ,[C_r]$ in $\Gamma$ such that 
$$
\overline{NE}_1(X)= \NEe + \sum_{i=1}^{r} \r_{\geq 0}[C_i], \text{ and}
$$
\item $\overline{NE}_1(X)= \NE + \sum_{[C]\in \Gamma} \r_{\geq 0}[C]$.
\end{enumerate}
\end{thm*}

We denote by $\overline{NM}_1(X)$ the closed convex cone in $N_1(X)$ generated by classes
of curves on $X$ with nonnegative intersection with all effective divisors on $X$.
By \cite[Theorem 2.2]{BDPP_peff_cone}, $\overline{NM}_1(X)$ 
is the closed convex cone in $N_1(X)$ generated by  classes of the 
so called \emph{strongly movable curves}.
Strongly movable curves are images  in $X$ of curves obtained as complete intersections of
suitable very ample divisors on birational modifications of $X$.
In \cite[Theorem 2.2]{BDPP_peff_cone} $X$ is assumed to be smooth, but one can easily deduce
the general case from the smooth case by a standard blowup argument.

The aim of this paper is to give a description of the cone $\overline{NM}_1(X)$ from the point of view of 
the Minimal Model Program. 
Let $(X,\Delta)$ be a $\q$-factorial klt  pair, and $\varphi:X\map X'$ a birational map  obtained by a Minimal Model Program.
We define the \emph{numerical pullback of curves} $\varphi_{num}^*:N_1(X')\to N_1(X)$
as the dual linear map of $\varphi_*:N^1(X)\to N^1(X')$.
Here $N^1(X)$ denotes  the $\r$-vector space $NS(X)\otimes \r$, which is identified with the dual vector space of $N_1(X)$.

Our main theorem is the following.

\begin{thm}  \label{generalization}
Let $(X,\Delta)$ be a $\q$-factorial klt  pair.
Let $\Sigma\subset N_1(X)$ be the set of classes of curves on $X$
that are  numerical pullbacks of curves lying on general fibers of
Mori Fiber Spaces obtained from $X$ by a $(K_X+\Delta)$-Minimal Model Program with scaling.
Then 
$$
\NE + \nm = \NE +\overline{\sum_{[C]\in \Sigma}\r_{\geq 0}[C]}.
$$
\end{thm}

\begin{cor} \label{log_fano}
Let $(X,\Delta)$ be a $\q$-factorial klt  pair such that $-(K_X+\Delta)$ is ample.
Then there are finitely many curves  $C_1,\dots ,C_r\subset X$ such that 
 \begin{enumerate}
\item[$\bullet$] each $[C_i]$ is the numerical pullback of a curve lying on a general fiber of
a Mori Fiber Space obtained from $X$ by a $(K_X+\Delta)$-Minimal Model Program with scaling,  and
\item[$\bullet$] $\nm = \sum_{i=1}^{r}\r_{\geq 0}[C_i]$.
\end{enumerate}
\end{cor}

For a $\q$-factorial klt pair $(X,\Delta)$, and an ample $\r$-divisor $H$ on $X$,  
the effective threshold $\sigma(H)$ of $H$ is the supremum of the set of real numbers $t$ for which
$H+t(K_X+\Delta)$ is pseudo-effective. 
In the course of the proof of the main theorem, we show in Theorem~\ref{thm:ample}
that $\sigma(H)$ is a rational number when both $\Delta$ and $H$ have rational coefficients.
When $\Delta=\emptyset$ and $X$ is a terminal threefold, this was proved by Batyrev in \cite{batyrev}.
When $(X,\Delta)$ is a $3$-dimensional $\q$-factorial  klt pair, 
this was proved by Fujita in  \cite{fujita_kodaira_energy}.

In \cite{batyrev}, Batyrev proved Theorem~\ref{generalization} in the special case when 
$\Delta=0$ and $X$ is a terminal threefold.
Moreover he claimed  the following.

\begin{thm} \label{main_thm}
Let $X$ be a projective $\q$-factorial terminal threefold.
Let $\Sigma\subset N_1(X)$ be as defined in Theorem~\ref{generalization}, with $\Delta=0$.
\begin{enumerate}
\item If $A$ is an ample $\r$-divisor on $X$, then  the set of rays 
$\big\{\r_{\geq 0}[C]\ \big|\ [C]\in \Sigma \text{ and } (K_X+ A)\cdot C<0\big\}$ is finite.
In particular, there are finitely many classes  $[C_1],\dots ,[C_r] \in \Sigma$ such that
$$
\NEb + \nm = \Big(\NEb + \nm \Big)_{_{(K_X+A)\geq 0}}+\sum_{i=1}^{r}\r_{\geq 0}[C_i].
$$
\item $\NEb + \nm = \NEb +\sum_{[C]\in \Sigma} \r_{\geq 0}[C]$.
\end{enumerate}
\end{thm}

Batyrev's  proof, however, contains a gap.
Namely, in \cite{batyrev} it is assumed and used that 
numerical pullbacks of curves on general fibers of Mori Fiber Spaces obtained from $X$ 
by a $(K_X+\Delta)$-Minimal Model Program are always integral. 
This is not  true though, as Example~\ref{counter_example} illustrates.
To the best of our knowledge, this is the first
complete proof of this result. In order
to complete Batyrev's proof, we rely on boundedness of terminal
Fano threefolds. For this reason, as it stands, this proof does
not generalize to the log terminal case.

Strictly speaking, Theorems \ref{generalization} and \ref{main_thm} are structure theorems for the cone
$\NE + \nm$ rather than the cone $\nm$. 
In fact, the analogous statements do not hold in general for the cone $\nm$,
as the next example illustrates.

\begin{exmp}
Let $C\subset \p^2$ be a smooth cubic curve, $P\in C$ an inflexion point, and
$L\subset \p^2$ the tangent line to $C$ at $P$. Consider the pencil of cubic curves on $\p^2$ generated
by $C$ and $3L$. All base points of this pencil lie over $P$. Let $X$ be the surface obtained
from $\p^2$ by blowing up 8 of the base points.
Let $\tilde C$ and $\tilde L\subset X$ be the birational transforms of $C$ and $L$ respectively.
Then $K_X\cdot\tilde C<0$.
For each   $1\leq i\leq 7$, let $E_i\subset X$ be the the birational transform of the exceptional 
curve of the  $i$-th blowup.
Notice that $\tilde L\cdot \tilde C=E_1\cdot \tilde C= \cdots =E_7\cdot \tilde C=0$.
Since $\dim N^1(X)=9$ and $[\tilde L], [E_1], \dots, [E_7]$
are linearly independent in $N^1(X)$, it follows that 
$\r[\tilde C]=\big\{z\in N_1(X)\ \big|\ \tilde L\cdot z=E_1\cdot z= \cdots =E_7\cdot z=0\big\}$.
Notice also that, $\forall z\in \nm$, $\tilde L\cdot z\geq 0$ and $E_i\cdot z\geq 0$ for $1\leq i\leq 7$.
Thus  $\r_{\geq0}[\tilde C]$ is an extremal ray of the cone $\nm$.

On the other hand, $\r_{\geq0}[\tilde C]$ is not generated by 
the numerical pullback of any curve on a fiber of a Mori Fiber Space $f':X'\to Y'$ obtained 
from $X$ by a Minimal Model Program.
Indeed, assume otherwise and let $C'$ be the image of $\tilde C$ in $X'$.
Then $C'$ lies on a fiber of $f'$. Since $X$ is rational and $C$ has genus $1$, this implies that $X'\cong \p^2$.
Moreover, $\tilde L, E_1, \dots, E_7$ are precisely the curves on $X$ contracted by the morphism $X\to X'$.
But this is impossible since none of these is a $-1$-curve.
\end{exmp}

We also refer to \cite[Example 5.4]{lehmann_nef_cone} for an example of a smooth threefold $X$ and ample $\r$-divisor $A$
for which 
$
\nm \neq \nm_{_{(K_X+A)\geq 0}}+\sum \r_{\geq 0}[C_i]
$
for any countable collection of integral classes of curves $\{[C_i]\}$.

This paper is organized as follows. 
In section~\ref{background} we define the cones of curves and divisors that will be
of interest to us. We also define nef and effective thresholds.
In section~\ref{section_mmp}  we explain the Minimal Model Program with scaling established in \cite{bchm}.
In section~\ref{mistake} we define numerical pullback of curves and discuss examples.
In section~\ref{sigma_is_rational}, using the Minimal Model Program with scaling, 
we prove Theorem~\ref{generalization} and Corollary~\ref{log_fano}.
In the course of the proof, we also prove rationality of  the effective threshold of ample 
$\q$-divisors on $X$ when $\Delta$ has rational coefficients.
In section~\ref{end} we prove Theorem~\ref{main_thm}.

We work over some fixed algebraically closed field of characteristic zero,
and follow the terminology in \cite{kollar_mori} and \cite{kamama}.
Curves are always assumed to be projective.

\medskip
{\bf Acknowledgments.}  
This paper grew from the older notes \cite{after_batyrev}.
Using recent results from \cite{bchm}, the arguments presented there can be extended 
to arbitrary $\q$-factorial klt  pairs.
I am deeply grateful to St\'ephane Druel for very fruitful discussions, 
for suggesting several simplifications and improvements of the proofs,  
in particular pointing out Lemma~\ref{lambda->0} and Remark~\ref{exposed} to me.
I am greatly indebted to the referee for pointing out many inaccuracies and imprecisions in previous 
versions of these notes and suggesting many improvements. 
This paper has benefitted immensely from their  comments. 
I also thank J\'anos Koll\'ar for his comments and suggestions. 
After this work had been completed, I have learned that Brian Lehmann
has independently obtained Theorem~\ref{generalization} in \cite{lehmann_nef_cone}.


\section{Cones of curves and divisors and nef and effective thresholds}
\label{background}

In this section we define some important cones of curves and divisors.
We also define nef and effective thresholds.
We refer to \cite{kollar_mori} and references therein for details and proofs.

Let $X$ be a normal projective variety.

\begin{defn} \label{cones} 
Let $N^1(X)$ be   the $\r$-vector space $NS(X)\otimes \r$, where $NS(X)$ denotes the
N\'eron-Severi group of $X$.
Given an $\r$-linear combination $D$ of Cartier divisors on $X$, we denote by $[D]$
the corresponding element in $N^1(X)$.
A $1$-cycle on $X$ is a formal $\r$-linear combination of irreducible and reduced curves on $X$.
Two $1$-cycles $C$ and $C'$ are said to be numerically equivalent if $D\cdot C=D\cdot C'$
for every Cartier divisor $D$ on $X$.
Let $N_1(X)$ be the $\r$-vector space of $1$-cycles on $X$  modulo numerical equivalence.
By abuse of notation, we identify curves on $X$ with the $1$-cycles naturally associated to them.
Given a $1$-cycle $C$ on $X$, we denote by $[C]$ its class in $N_1(X)$.
Intersection number of Cartier divisors and curves
induces a perfect pairing between $N^1(X)$ and $N_1(X)$.
These two vector spaces are finitely generated. 
Their common dimension is denoted by $\rho(X)$, the Picard number of $X$.

The \emph{cone of curves} of $X$, $\overline{NE}_1(X)$, is
defined to be the closed convex cone in $N_1(X)$ generated by the classes
of irreducible curves on $X$. We denote its dual cone by $\nef(X)\subset N^1(X)$.
It is the closed convex cone generated by the classes of nef divisors on $X$.

Similarly, we define the \emph{cone of pseudo-effective divisors} of $X$,
$\eff(X)$, to be the closed convex cone in $N^1(X)$ generated by the classes of effective divisors on $X$. 
We denote the dual cone of $\eff(X)$ by $\overline{NM}_1(X)\subset
N_1(X)$. It is called the \emph{cone of nef curves} of $X$.
\end{defn}

\begin{defn} \label{q-divisors}
A \emph{$\q$-divisor} on $X$ is a $\q$-linear combination of prime Weil divisors on $X$. 
A $\q$-divisor $D$ on $X$ is said to be \emph{$\q$-Cartier} if some nonzero multiple of $D$
is a Cartier divisor. 
Two $\q$-divisors $D$ and $D'$ on $X$ are said to be \emph{$\q$-linearly equivalent} if
there exists an integer $m>0$ such that both $mD$ and $mD'$ are Cartier and 
$mD\sim mD'$. In this case we write $D\sim_{\q}D'$.

We say that  $X$ is \emph{$\q$-factorial} if every  $\q$-divisor on $X$
is $\q$-Cartier.
\end{defn}

\begin{defn} \label{big_divisors}
A $\q$-Cartier $\q$-divisor $D$ on $X$ is said to be \emph{big} if there exists a constant $c>0$ such that
$h^0(X,\ou_X(kD))\geq c\cdot k^{\dim(X)}$ for all $k$ sufficiently large and divisible.
Equivalently, $D$ is big if and only if the linear system $\big|\ou_X(kD)\big|$ defines a 
birational map for all $k$ sufficiently large and divisible.
\end{defn}

\begin{rem} \label{big}
Let $D$ be a $\q$-Cartier $\q$-divisor on $X$.
By Kleiman's ampleness criterion \cite{kleiman66},
$D$ is ample if and only $[D]\in \inter\big(\nef(X)\big)$.
Therefore $\nef(X)\subset \eff(X)$,
and $\overline{NM}_1(X)\subset \overline{NE}_1(X)$.

By Kodaira's lemma,  $D$ is big if and only if there are ample $\q$-Cartier $\q$-divisor $H$ 
and effective $\q$-Cartier $\q$-divisor $E$ such that $D\sim_{\q}H+E$. 
Thus $D$ is big if and only if  $[D]\in \inter\big(\eff(X)\big)$.
\end{rem}

\begin{defn} \label{r-divisors}
An \emph{$\r$-divisor} on $X$ is an $\r$-linear combination of prime Weil divisors on $X$. 
An $\r$-divisor on $X$ is said to be \emph{$\r$-Cartier} if it  is an $\r$-linear combination 
of Cartier divisors on $X$. 
Two $\r$-divisors $D$ and $D'$ on $X$ are said to be \emph{$\r$-linearly equivalent} if
there exist real numbers $a_1,\dots, a_k$ and  rational functions $f_1, \dots, f_k$ on $X$ such that 
$D-D'=\sum_{i=1}^k a_i \ div(f_i)$.
In this case we write $D\sim_{\r}D'$.

Let $D$ be an $\r$-Cartier $\r$-divisor on $X$.
We say that $D$ is \emph{ample} (respectively  \emph{big}, respectively \emph{pseudo-effective})
if $[D]\in \inter\big(\nef(X)\big)$ (respectively  $[D]\in \inter\big(\eff(X)\big)$, respectively $[D]\in \eff(X)$).
\end{defn}

\begin{defn}\label{extremal_rays}
An \emph{extremal face} $F$ of a cone $N\subset \r^n$ is a subcone of $N$
satisfying:
$$
u,v\in N \ \text{ and } \ u+v\in F \ \ \Rightarrow \ \ u,v\in F.
$$
A $1$-dimensional extremal face of $N$ is called an \emph{extremal
ray}.

Let $D:\r^n\to \r$ be a linear function. Set $N_{D\geq0}=\{z\in N\ |\ D(z)\geq 0\}$, 
and similarly for $N_{D=0}$, $N_{D\leq0}$, etc.
An extremal face $F\subset N$ such that $F\setminus\{0\} \subset N_{D<0}$
is called a \emph{$D$-negative extremal face}.
If $F\subset N_{D=0}$, then we say that \emph{$F$ is supported on $D$}.

When dealing with $N_1(X)$, by abuse of notation, we identify an element  $\alpha\in N^1(X)$ with 
the linear function $N_1(X)\to \r$ defined by $\alpha$ via intersection product.
\end{defn}

\begin{defn} \label{thresholds}
Let  $\Delta$ be an $\r$-divisor on $X$ such that $K_X+\Delta$ is $\r$-Cartier, and $H$
an $\r$-Cartier $\r$-divisor on $X$. We define the \emph{nef  threshold of $H$} by
$$
\tau(X,\Delta,H)= \sup\Big\{t\geq 0\ \Big| \ \big[H+t(K_X+\Delta)\big] \in \nef(X) \Big\},
$$
provided that the set on the right hand side is nonempty.
We define the \emph{effective threshold of $H$} by
$$
\sigma(X,\Delta,H)= \sup\Big\{t\geq 0\ \Big|\ \big[H+t(K_X+\Delta)\big] \in \eff(X)\Big\},
$$
provided that the set on the right hand side is nonempty.
When there is no ambiguity we write $\tau(H)$ and $\sigma(H)$ for
$\tau(X,\Delta,H)$ and $\sigma(X,\Delta,H)$ respectively.
\end{defn}

\begin{rem}\label{tau<sigma}
Let $\Delta$ be as in Definition~\ref{thresholds}.
Suppose that $\big[H+t(K_X+\Delta)\big] \in \nef(X)$ for some $t\geq 0$. 
Then, since $\nef(X)\subset \eff(X)$,  
$\tau(H)\leq \sigma(H)$.
\end{rem}

\begin{rem}
In \cite{fujita_kodaira_energy} Fujita defined the log Kodaira Energy
$\kappa\epsilon(X,\Delta,H)$. This is related to the effective
threshold $\sigma(X,\Delta,H)$ by the formula
$$
\sigma(X,\Delta,H)=-\frac{1}{\kappa\epsilon(X,\Delta,H)}.
$$
\end{rem}


\section{The Minimal Model Program with scaling}
\label{section_mmp}
In this section we gather some results about the Minimal Model Program - MMP for short - and fix notation.
We explain in detail the MMP with scaling established in \cite{bchm}, which 
will be essential in the proof of our main theorem.

\begin{defn}[{See \cite[section 2.3]{kollar_mori}}] \label{def_klt}
Let $X$ be a normal projective variety, and $\Delta=\sum a_iD_i$ an effective  $\r$-divisor on $X$,
where the $D_i$'s are distinct prime divisors. Suppose that $K_X+\Delta$ is $\r$-Cartier. 
Let $f:\tilde X\to X$ be a log resolution of the pair $(X,\Delta)$. 
This means that $\tilde X$ is a smooth projective 
variety, $f$ is a birational morphism whose exceptional locus is the union of prime divisors $E_i$'s, 
and the divisor $\sum E_i+f^{-1}_*\Delta$ has simple normal crossing 
support. 
There are uniquely defined real numbers $a(E_i,X,\Delta)$'s such that
$$
K_{\tilde X}+f^{-1}_*\Delta\sim_{\r} f^*(K_X+\Delta)+\sum_{E_i}a(E_i,X,\Delta)E_i.
$$
The $a(E_i,X,\Delta)$'s do not depend on the log resolution $f$,
but only on the valuations associated to the $E_i$'s.

We say that $(X,\Delta)$, or $K_X+\Delta$, is \emph{klt} if  $0<a_i<1$  and, for some  log resolution 
$f:\tilde X\to X$ of $(X,\Delta)$, $a(E_i,X,\Delta)>-1$ 
for every $f$-exceptional prime divisor $E_i$.
If this condition holds for some log resolution of $(X,\Delta)$, then it holds for every  
log resolution of $(X,\Delta)$.

When $X$ is a $\q$-factorial normal projective variety and $(X,\Delta)$ is klt, 
we say that $(X,\Delta)$ is  a \emph{$\q$-factorial klt pair}.
We stress that $\Delta$ is allowed to have real coefficients.
\end{defn}

\begin{thm}[{Rationality Theorem - \cite[Lemma 2.1]{keel_matsuki_mckernan_abu}}] \label{rationality}
Let $(X,\Delta)$ be a $\q$-factorial klt pair, where $\Delta$ is a $\q$-divisor, and 
let $H$ be a $\q$-divisor on $X$.
Suppose that $K_X+\Delta$ is not nef and that $H+t(K_X+\Delta)$ is nef for some $t\in \q_{\geq 0}$. 
Then $\tau(X,\Delta,H)\in \q$. 
\end{thm}

\begin{lemma}[{\cite[Lemma 2.6]{birkar07}}]\label{exists_extremal_ray}
Let $(X,\Delta)$ be a $\q$-factorial klt pair.
Suppose that $K_X+\Delta$ is not nef. Let $H$ be an effective
$\r$-divisor on $X$ such that  $K_X+\Delta+H$ is nef and klt. 
Then there is a $(K_X+\Delta)$-negative extremal ray $R$ of the cone $\nec(X)$ 
supported on $H+\tau(H)(K_X+\Delta)$. 
\end{lemma}

Let $(X,\Delta)$ be a $\q$-factorial klt pair and $F\subset \nec(X)$ a 
$(K_X+\Delta)$-negative extremal face.
Then there is a $\q$-divisor $\Delta'$ on $X$ such that $(X,\Delta')$ is klt  and $F$ is 
$(K_X+\Delta')$-negative.
So one applies \cite[Theorem 3.7.3]{kollar_mori} to get the following.

\begin{thm}[{Contraction Theorem}]\label{contraction_thm}
Let $(X,\Delta)$ be a $\q$-factorial klt pair.
Let $F$ be a $(K_X+\Delta)$-negative extremal face of  the cone $\nec(X)$.
Then there exists a unique morphism $f:X\to Y$ onto a normal projective variety such that 
$f_*\ou_{X}=\ou_{Y}$, and,
for any curve $C\subset X$, $\varphi(C)$ is a point if and only if $[C]\in F$.  
\end{thm}

\begin{defn}
Under the assumptions and notation of Theorem~\ref{contraction_thm}, 
we say that $f:X\to Y$ is \emph{the contraction associated to} $F$.
\end{defn}

\begin{say}[Properties of contractions associated to extremal rays] \label{rays}
Let $(X,\Delta)$ be a $\q$-factorial klt pair.
Let $R$ be a $(K_X+\Delta)$-negative extremal ray of the cone $\nec(X)$,
and let $f:X\to Y$ be the contraction associated to $R$.
Then there is a $\q$-divisor $\Delta'$ on $X$ such that $(X,\Delta')$ is klt  and $R$ is $(K_X+\Delta')$-negative.
By  \cite[Theorem 3.7]{kollar_mori}, there is an exact sequence
$0\to \pic(Y) \xrightarrow{f^*}\pic(X)\to\z$, where the last map is given by intersecting 
with a curve $C\subset X$ such that $R=\r_{\geq 0}[C]$. 
Thus $\rho(X)=\rho(Y)+1$.
By tensoring this sequence  with $\r$, we see that, given an
$\r$-divisor $D$ on $X$ such that $D\cdot R=0$, 
there is an $\r$-Cartier $\r$-divisor $D_Y$ on $Y$ such that $D\sim_{\r} f^*D_Y$.

By \cite[Proposition 2.5]{kollar_mori}, one of the following  situations occurs.
\begin{enumerate}
  \item $\dim(Y)<\dim(X)$. 
     We call such $f$ a \emph{Mori Fiber Space}, and we say that $R$ is a \emph{fiber type extremal ray}.
  \item The morphism $f$ is birational and the exceptional locus of $f$ consists of a prime divisor on $X$.
    In this case, $(Y,f_*\Delta)$ is a $\q$-factorial klt pair.
    We call such $f$ a \emph{divisorial contraction}.
  \item The morphism $f$ is birational and the exceptional locus of $f$ 
    has codimension at least $2$ in $X$.
    We call such $f$ a \emph{small contraction}.
\end{enumerate}
\end{say}

The following result was established in \cite[Corollary 1.4.1]{bchm}.

\begin{defn-thm}[{Existence of flips}] \label{existence_of_flips}
Let $(X,\Delta)$ be a $\q$-factorial klt pair.
Let $f:X\to Y$ be a small contraction associated to a $(K_X+\Delta)$-negative extremal ray 
$R\subset \nec(X)$.
Then there exist unique birational map $\varphi:X\map X^+$ and morphism $f^+:X^+\to Y$ 
such that the following diagram commutes
\[
\xymatrix{
X \ar@{>}[dr]_{f} \ar@{.>}[rr]^{\varphi} & & X^+ , \ar@{>}[dl]^{f^+} \\
& Y
}
\]
$(X^+,{\varphi}_*\Delta)$ is a $\q$-factorial klt pair,
the exceptional locus of $f^+$ has codimension at least $2$ in $X^+$,
and $(K_{X^+}+{\varphi}_*\Delta)\cdot C>0$ for every curve $C\subset X^+$ contracted by $f^+$.
We say that the map $\varphi:X\map X^+$ is \emph{the flip associated to} $R$.
\end{defn-thm}

For a general description of the MMP, we refer the reader to \cite[sections 2.1 and 2.2]{kollar_mori}.
Here we explain how to use an ample $\r$-divisor $H$ on $X$
to drive the MMP.  
This is what is called   the $(K_X+\Delta)$-MMP with scaling of $H$.

\begin{say} [{The MMP with scaling}] \label{MMP_with_scaling}
Let $(X,\Delta)$ be a $\q$-factorial klt pair, and let $H$ be an effective $\r$-divisor on $X$
such that $K_X+\Delta+H$ is nef and klt. 
We will define inductively (possibly finite) sequences of projective varieties $X_i$, birational maps 
$\varphi_i:X_i\map X_{i+1}$,  $\r$-divisors $\Delta_i$ and $H_i$ on $X_i$, 
$(K_{X_i}+\Delta_i)$-negative extremal rays $R_i\subset \nec(X_i)$, 
and  real numbers $0\leq \lambda_i\leq 1$.
\medskip
 
\noindent {\bf Step 0.} We set $X_0=X$, $\Delta_0=\Delta$, and $H_0=H$.
By hypothesis, $K_{X_0}+\Delta_0 + H_0$ is nef and klt.
We move to Step 1 with $n=0$
\medskip

\noindent {\bf Step 1.} Suppose we have constructed $X_n$, $\Delta_n$ and $H_n$ such that
 $(X_n,\Delta_n)$ is a $\q$-factorial klt pair, and
$K_{X_n}+\Delta_n + \lambda H_n$ is nef and klt for some $0\leq \lambda \leq 1$.
Set $\lambda_n=\inf\Big\{\lambda\geq 0 \ \Big|\ \big[K_{X_n}+\Delta_n + \lambda H_n\big]\in \nef(X_n)\big\}$.
Then $0\leq \lambda_n\leq 1$, and $K_{X_n}+\Delta_n + \lambda_n H_n$ is nef and klt. 
We move to Step 2. 
\medskip

\noindent {\bf Step 2.} We check whether  $K_{X_n}+\Delta_n$ is nef. 

If $K_{X_n}+\Delta_n$ is nef, then we stop and the sequence $\big\{X_i\big\}$ ends with $X_n$.

If $K_{X_n}+\Delta_n$ is not nef, then, by Lemma~\ref{exists_extremal_ray},  there exists at least one 
$(K_{X_n}+\Delta_n)$-negative extremal ray $R\subset \nec(X_n)$ such that 
$(K_{X_n}+\Delta_n + \lambda_n H_n)\cdot R=0$.
We choose one such extremal ray $R_n$, and let $f:X_n\to Y$ be the contraction associated to $R_n$.
We move to Step 3. 
\medskip

\noindent {\bf Step 3.} 
We check which of the three possibilities described in \ref{rays} occurs. 
\smallskip

 {\bf (1)} If $\dim Y<\dim X_n$, then we stop and the sequence $\big\{X_i\big\}$ ends with $X_n$. 
\smallskip

 {\bf (2)}  If $f:X_n\to Y$ is a divisorial contraction, then
we set $X_{n+1}=Y$, $\varphi_n=f$, $\Delta_{n+1}=f_*\Delta_n$, and $H_{n+1}=f_{*}H_n$.
So $(X_{n+1},\Delta_{n+1})$ is a $\q$-factorial klt pair.
Since $(K_{X_n}+\Delta_n + \lambda_n H_n)\cdot R_n=0$, by \ref{rays},
\begin{equation} \label{eq_discrep}
K_{X_{n}}+\Delta_{n} + \lambda_n H_{n}\sim_{\r}
(\varphi_n)^*\big(K_{X_{n+1}}+\Delta_{n+1} + \lambda_n H_{n+1}\big).
\end{equation}
This  implies that $K_{X_{n+1}}+\Delta_{n+1} + \lambda_n H_{n+1}$  is nef and klt
since so is $K_{X_{n}}+\Delta_{n} + \lambda_n H_{n}$.
We go back to Step 1 replacing $n$ with $n+1$.
\smallskip

{\bf (3)} If $f:X_n\to Y$ is a small contraction, and $\varphi:X_n\map X_n^+$ is the associated flip, then 
 we set $X_{n+1}=X_n^+$, $\varphi_n=\varphi$, $\Delta_{n+1}=\varphi_*\Delta_n$, and $H_{n+1}=\varphi_*H_n$.
So $(X_{n+1},\Delta_{n+1})$ is a $\q$-factorial klt pair.
Consider the flip diagram:
\[
\xymatrix{
X_n \ar@{>}[dr]_{f} \ar@{.>}[rr]^{\varphi_n} & & X_{n+1} \ar@{>}[dl]^{f^+} \\
& Y
}
\]
Since $(K_{X_n}+\Delta_n + \lambda_n H_n)\cdot R_n=0$, by \ref{rays}, 
there exists an $\r$-Cartier $\r$-divisor $D_Y$ on $Y$ such that
$K_{X_n}+\Delta_n + \lambda_n H_n\sim_{\r} f^*D_Y$. 
Then $K_{X_{n+1}}+\Delta_{n+1} + \lambda_n H_{n+1}\sim_{\r}(f^+)^*D_Y$.
By hypothesis $K_{X_n}+\Delta_n + \lambda_n H_n$ is nef and klt.
Thus $D_Y$ is nef and so is $K_{X_{n+1}}+\Delta_{n+1} + \lambda_n H_{n+1}$.
By looking at a common log resolution of $\big(X_n,\Delta_n+\lambda_n H_n\big)$ and 
$\big(X_{n+1},\Delta_{n+1} + \lambda_n H_{n+1}\big)$, 
we conclude that $K_{X_{n+1}}+\Delta_{n+1} + \lambda_n H_{n+1}$
is also klt.
We go back to Step 1 replacing $n$ with $n+1$.
\medskip

By abuse of notation, we denote the set of six sequences constructed above by 
$\big\{(X_i, \varphi_i, \Delta_i, H_i, \lambda_i, R_i)\big\}$.
They satisfy the following conditions, whenever the corresponding objects are defined.
\begin{enumerate}
	\item[$\bullet$] $\lambda_i=\inf\Big\{\lambda\geq 0 \ \Big|\ \big[K_{X_i}+\Delta_i + \lambda H_i\big]\in \nef(X_i)\Big\}$, and  $\big\{\lambda_i\big\}$ is a decreasing sequence of  real numbers
		between $0$ and $1$. If $\Delta$ and $H$ are $\q$-divisors, then $\lambda_i\in \q$ by Theorem~\ref{rationality}.
	\item[$\bullet$] $R_i$ is a $(K_{X_i}+\Delta_i)$-negative extremal ray of $\nec(X_i)$ such that  $(K_{X_i}+\Delta_i + \lambda_i H_i)\cdot R_i=0$.
	\item[$\bullet$] $\varphi_i:X_i\map X_{i+1}$  is the divisorial contraction or flip associated to  $R_i$, $\Delta_{i+1}=(\varphi_{i})_{*}\Delta_i$, and $H_{i+1}=(\varphi_{i})_{*}H_i$.
		If $H$ is big, then so is $H_i$.
\end{enumerate}

We say that  $\big\{(X_i, \varphi_i, \Delta_i, H_i, \lambda_i, R_i)\big\}$ is a 
\emph{$(K_X+\Delta)$-MMP with scaling of $H$}. 
If  the sequence $\big\{X_i\big\}$ ends with $X_n$ for some $n\geq 0$, then 
\begin{enumerate}
  \item[$\bullet$] either $K_{X_n}+\Delta_n$ is nef, or
  \item[$\bullet$] the extremal ray $R_n\subset \nec(X_n)$ is a fiber type extremal ray. 
In this case we say that the MMP \emph{terminates with the Mori Fiber Space}
  $f:X_n\to Y$ associated to $R_n$.
\end{enumerate}
\end{say}

By \cite[Corollary 1.3.2]{bchm}, if $K_X+\Delta\not \in \eff(X)$, then some
$(K_X+\Delta)$-MMP with scaling terminates with a Mori Fiber Space. 
In fact, the proof of \cite[Corollary 1.3.2]{bchm} gives the following stronger result.
Since it is not explicitly stated in \cite{bchm}, for the reader's convenience, we provide a proof below.

\begin{thm}\label{mmp_w_scaling_terminates}
Let $(X,\Delta)$ be a $\q$-factorial klt pair such that $K_X+\Delta\not \in \eff(X)$. 
Let $H$ be an effective ample $\r$-divisor  on $X$ such that $K_X+\Delta+H$ is nef and klt.
Then any $(K_X+\Delta)$-MMP with scaling of $H$ terminates with a Mori Fiber Space.
\end{thm}

The proof of Theorem~\ref{mmp_w_scaling_terminates} uses the finiteness of models established in \cite{bchm}.
In order to state the finiteness result, we need to define nef models. These are called 
weak log canonical models in \cite[Definition 3.6.6]{bchm}.

\begin{defn}\label{nef_models}
Let $(X,\Delta)$ be a $\q$-factorial klt pair.
A \emph{nef model} for $(X,\Delta)$ is a  $\q$-factorial klt pair $(X',\Delta')$,
together with a birational map $\varphi':X\map X'$ such that 
\begin{enumerate}
	\item $(\varphi')^{-1}$ does not contract any divisor,
	\item $\Delta'=(\varphi')_*\Delta$,
	\item $K_{X'}+\Delta'$ is nef, and
	\item for some, and hence all, common log resolution $p:W\to X$, $q:W\to X'$ of $(X,\Delta)$ and $(X',\Delta')$,
		$p^*(K_{X}+\Delta)=q^*(K_{X'}+\Delta')+E$, where $E$ is effective and $q$-exceptional.
\end{enumerate}
We say that two nef models $\big((X',\Delta'), \varphi'\big)$ and $\big((X'',\Delta''), \varphi''\big)$ 
for $(X,\Delta)$
are isomorphic if there is an isomorphism $\psi:X'\to X''$ such that $\varphi''=\psi\circ \varphi'$.
When the birational map $\varphi':X\map X'$ is understood from the context, we simply say that 
$(X',\Delta')$ is a nef model for $(X,\Delta)$.
\end{defn}

\begin{rem}\label{nef_models_under_sim}
Let $(X,\Delta)$ be a $\q$-factorial klt pair, 
and let $\tilde \Delta$ 
be another effective  $\r$-divisor on $X$ such that $(X,\tilde \Delta)$ is klt.
Suppose that $\Delta\sim_{\r} \tilde \Delta$ and let $\varphi:X\map X'$ be a birational map. Then 
$(X',\varphi_*\Delta)$ is a nef model for $(X,\Delta)$ if and only if $(X',\varphi_*\tilde \Delta)$ 
is a nef model for $(X,\tilde \Delta)$.
\end{rem}

\begin{rem}\label{scaling_yields_nef_models}
Let $(X,\Delta)$ be a $\q$-factorial klt pair. Let $H$ be an effective $\r$-divisor  on $X$
such that $K_X+\Delta+H$ is nef and klt.
Let $\big\{(X_i, \varphi_i, \Delta_i, H_i, \lambda_i, R_i)\big\}$ be a $(K_X+\Delta)$-MMP with scaling of 
$H$ as described in \ref{MMP_with_scaling}.
Then, for every $i\geq 0$ for which $X_i$ is defined, $(X_i,\Delta_i+\lambda_iH_i)$ is a nef model for 
$(X,\Delta+\lambda_iH)$.
Indeed, conditions (1)--(3) from Definition~\ref{nef_models} are straightforward.
To prove condition (4) we will show by induction on $i$ that, for some common log resolution 
$p:W\to X$, $q_i:W\to X_i$ of $(X,\Delta+H)$ and $(X_i,\Delta_i+H_i)$, and for every $0\leq t\leq \lambda_i$, 
\begin{equation} \label{(*)}
p^*(K_{X}+\Delta+tH)=q_i^*(K_{X_i}+\Delta_i+tH_i)+E_i(t),
\end{equation}
where $E_i(t)$ is effective and $q_i$-exceptional.

Suppose that \eqref{(*)} holds and that $\lambda_i>0$. 
By replacing $W$ with a further log resolution, we may assume that there is a log resolution 
$q_{i+1}:W\to X_{i+1}$  of  $(X_{i+1},\Delta_{i+1}+H_{i+1})$. 
The arguments in Step 3 of \ref{MMP_with_scaling} show that 
\begin{equation} \label{(**)}
q_i^*(K_{X_i}+\Delta_i+\lambda_iH_i)=q_{i+1}^*(K_{X_{i+1}}+\Delta_{i+1}+\lambda_iH_{i+1}).
\end{equation}
Since $\varphi_i$ is a step in the $(K_X+\Delta)$-MMP, 
there is an effective $q_{i+1}$-exceptional divisor $E_{i+1}$ on $W$
such that $q_i^*(K_{X_i}+\Delta_i)=q_{i+1}^*(K_{X_{i+1}}+\Delta_{i+1})+E_{i+1}$.
Together with \eqref{(**)} this implies that, for every $t\geq 0$, 
$$
q_i^*(K_{X_i}+\Delta_i+tH_i)=q_{i+1}^*(K_{X_{i+1}}+\Delta_{i+1}+tH_{i+1})+\frac{\lambda_i-t}{\lambda_i}E_{i+1}.
$$
Together with \eqref{(*)} this implies that, if  $0\leq t\leq \lambda_i$, 
then \eqref{(*)} holds with $i$ replaced with $i+1$. 
\end{rem}

The next finiteness result  is a special case of \cite[Theorem E]{bchm}.

\begin{thm}\label{finiteness}
Let $(X,\Delta)$ be a $\q$-factorial klt pair,  $H$ an effective ample $\r$-divisor  on $X$, 
and $\varepsilon>0$. Set
$$
\mathfrak{C}\ =\ \big\{ \tilde\Delta = \Delta +tH\ \big|\ t\geq \varepsilon \text{ and } (X,\tilde\Delta) \text{ is klt}  \big\}.
$$
Then the set of isomorphism classes of nef models  for pairs $(X,\tilde\Delta)$ 
such that $\tilde\Delta\in \mathfrak{C}$ is finite. 
\end{thm}

\begin{proof}
Write $\Delta=\sum c_i\Delta_i$ and $H=\sum a_iH_i$, where the $\Delta_i$'s 
are prime Weil divisors, the $H_i$'s are effective ample $\q$-divisors, and the $c_i$'s
and $a_i$'s are positive real numbers. 
Set $t_0= \sup\big\{t\geq 0\ \big| \ (X,\Delta +tH) \text{ is klt}  \big\}$.
We may assume that $t_0\geq \varepsilon$.
Fix rational numbers $\tilde a_i$ such that $0< \tilde a_i <\varepsilon a_i$.
Since the $H_i$'s are ample, we can choose ample $\q$-divisors $A_i\sim_{\q} H_i$ such that 
no $A_i$ has common components with any $\Delta_j$ or $H_j$, and such that 
$\big(X, \Delta + \sum \tilde a_iA_i + \sum (ta_i-\tilde a_i)H_i \big)$ is klt for every $\varepsilon\leq t<t_0$.
Set $A= \tilde a_iA_i$ and consider the finite dimensional affine subspace 
$V=\big\{ A+\sum t_i\Delta_i +\sum s_jH_j \ \big|\  t_i,s_j\in \r  \big\}$ of the 
real vector space  of Weil $\r$-divisors on $X$.
By \cite[Theorem E]{bchm}, the set of isomorphism classes of nef models for pairs $(X,\Delta')$ 
such that $\Delta'\in \cL_A$ is finite. 
Now notice that if $\tilde\Delta = \Delta +tH\in \mathfrak{C}$, then 
$\Delta' = A + \sum c_i\Delta_i + \sum (ta_i-\tilde a_i)H_i \in \cL_A$ and
$\Delta'\sim_{\r}\tilde\Delta$.
 The result then follows from  Remark~\ref{nef_models_under_sim}.
\end{proof}

The following is an immediate consequence of Theorem~\ref{finiteness} and 
Remark~\ref{scaling_yields_nef_models}.

\begin{lemma}\label{lambda->0}
Let $(X,\Delta)$ be a $\q$-factorial klt pair,  and let $H$ be an effective ample $\r$-divisor  on $X$
such that $K_X+\Delta+H$ is nef and klt.
Let $\big\{(X_i, \varphi_i, \Delta_i, H_i, \lambda_i, R_i)\big\}$ be a $(K_X+\Delta)$-MMP 
with scaling of $H$.
If this sequence is infinite, then $\lim_{i\to \infty}\lambda_i=0$. 
\end{lemma}

\begin{proof}[{Proof of Theorem~\ref{mmp_w_scaling_terminates}}]
Let $\big\{(X_i, \varphi_i, \Delta_i, H_i, \lambda_i, R_i)\big\}$ be a $(K_X+\Delta)$-MMP 
with scaling of $H$.
Suppose that this sequence is infinite.
Then  $\lim_{i\to \infty}\lambda_i=0$ by Lemma~\ref{lambda->0}.
On the other hand, it follows from  \eqref{(*)} that $K_X+\Delta+\lambda_i H$ 
is pseudo-effective for every $i$.
Thus $K_X+\Delta = \lim_{i\to \infty} K_X+\Delta+\lambda_i H$ is pseudo-effective, 
which is a contradiction.
So we conclude that the sequence $\big\{X_i\big\}$ ends with $X_n$ for some $n\geq 0$.
Moreover $\lambda_n>0$, which implies that $K_{X_n}+\Delta_n$ is not nef. 
Thus $R_n\subset \nec(X_n)$ is a fiber type extremal ray.
\end{proof}


\section{Numerical Pullback of Curves}
\label{mistake}

In this section we define numerical pullback of curves.

Let $\varphi:X\map Z$ be a birational map between $\q$-factorial
projective varieties, and assume that $\varphi^{-1}$ does not contract any divisor. 
Taking pullback of divisors on $Z$ defines an injective linear map
$\varphi^{*}:N^1(Z)\to N^1(X)$.
Taking pushforward of divisors on $X$ defines a surjective linear map
$\varphi_*:N^1(X)\to N^1(Z)$.
The composition $\varphi_*\circ \varphi^{*}$ is the identity on $N^1(Z)$.

\begin{defn}
We define the \emph{numerical pullback of curves} $\varphi_{num}^*:N_1(Z)\to N_1(X)$
as the dual linear map of $\varphi_*:N^1(X)\to N^1(Z)$. It is the unique 
injective linear map having the following properties.
\begin{enumerate}
\item[$\bullet$]  If $z\in N^1(Z)$ and $l\in N_1(Z)$, then
$\varphi^{*}(z)\cdot \varphi_{num}^*(l)=z\cdot l$.
\item[$\bullet$] If $\beta\in \ker \varphi_*$ and $m\in \im \varphi_{num}^*$, then
$\beta \cdot m =0$.
\end{enumerate}
\end{defn}

This definition allows us to take pullbacks of curves  
on $Z$ that are not the image of any curve on $X$, as in Example~\ref{flop} below. 
One must be careful though, as 
numerical pullbacks of integral (respectively effective) curves need not be integral 
(respectively effective), as the next two examples illustrate.

\begin{exmp}\label{counter_example}
Let $Y\subset \p^6$ be the cone over the Veronese surface in $\p^5$.
Then $Y$ is a $\q$-factorial terminal Fano threefold
of Picard number $1$.
Let $\pi:X\to Y$ be the blowup of the vertex of the cone,
and let $E\subset X$ be the exceptional divisor of $\pi$. Then $X$ is a smooth threefold and $E\cong \p^2$.
Let $l\subset Y$ be a ruling of the cone.
A simple computation shows that
$\pi_{num}^*\big([l]\big)=[\tilde l]+ \frac{1}{2}[e]$,
where $\tilde l$ is the strict transform of $l$ and
$e$ is a curve on $E$ corresponding to a line on $\p^2$
under the isomorphism $E\cong \p^2$.
\end{exmp}

\begin{exmp}\label{flop}
Let $Y\subset \p^4$ be the cone over a smooth quadric surface $Q\subset \p^3$. 
Let $L_1$ and $L_2$ denote $2$-planes on $Y$ coming from the two distinct rulings of $Q$.
For $i=1,2$,
let $\pi_i:X_i\to Y$ be the blowup of $Y$ along $L_i$. 
Then $X_i$ is a smooth threefold, and the exceptional locus of $\pi_i$ is a smooth 
rational curve $l_i$, which is contracted to the vertex of the cone.
Let $\varphi:X_1\map X_2$ be the flop $\pi_2^{-1}\circ \pi_1$. 
Then $\varphi_{num}^*\big([l_2]\big)=-[l_1]$, which is not effective.
\end{exmp}


\section{The structure of the cone of nef curves}
\label{structure_cone_pseff}
\label{sigma_is_rational}

In this section, using the MMP with scaling and ideas of \cite{batyrev}, we prove Theorem~\ref{generalization}
and Corollary~\ref{log_fano}.
Let  $(X,\Delta)$ be a $\q$-factorial klt pair.
If $K_X+\Delta \in \eff(X)$, then Theorem~\ref{generalization} is obviously true.
So we assume that $K_X+\Delta\not \in \eff(X)$. 
Let $H$ be an effective and ample $\r$-divisor  on $X$ such that $K_X+\Delta+H$ is nef and klt.
By Theorem~\ref{mmp_w_scaling_terminates},  any $(K_X+\Delta)$-MMP with scaling of $H$
terminates with a Mori Fiber Space $f:X'\to Y'$.
We will show in Theorem~\ref{thm:ample}  that
the numerical pullback to $X$ of any curve on a fiber of $f$ is supported on $H+\sigma(H)(K_X+\Delta)$.
This will be essential in the proof of Theorem~\ref{generalization}.

We start by investigating what happens with the effective thresholds when we run a $(K_X+\Delta)$-MMP with scaling of $H$.

\begin{lemma} \label{process}
Let $(X,\Delta)$ be a $\q$-factorial klt pair, and let $H$ be an effective  $\r$-divisor  on $X$ such that $K_X+\Delta+H$ is nef and klt.
Let $\big\{(X_i, \varphi_i, \Delta_i, H_i, \lambda_i, R_i)\big\}$  be a $(K_X+\Delta)$-MMP with scaling of $H$
as described in \ref{MMP_with_scaling}.
Then $\sigma(X_i,\Delta_i,H_i)=\sigma(X,\Delta,H)$ for every $i\geq 0$ for which $X_i$ is defined.
\end{lemma}

\begin{proof}
We will prove the result  by induction on $i$. 
So suppose $\sigma_{_n}:=\sigma(H_n)=\sigma(H)$ and consider $\varphi_n:X_n\map X_{n+1}$ in the sequence.

For any birational map $\varphi:X\map X'$ between $\q$-factorial projective varieties whose 
inverse does not contract any divisor, and for every pseudo-effective $\r$-divisor $D$ on $X$,
$\varphi_*D$ is also pseudo-effective. Thus $\sigma_{_{n+1}}\geq \sigma_{_n}$.
Moreover, if $\varphi_n$ is a flip, then $\sigma_{_{n+1}}=\sigma_{_n}$.
So we may assume that $\varphi_n:X_n\to X_{n+1}$ is a divisorial contraction. 
Let $E$ be the exceptional divisor of $\varphi_n$. 
Since $\varphi_n$ is a step in the $(K_X+\Delta)$-MMP, 
$K_{X_{n}}+\Delta_{n}=(\varphi_n)^*\big(K_{X_{n+1}}+\Delta_{n+1}\big)+aE$ for some $a>0$.
Together with \eqref{eq_discrep} and the fact that $\tau_{_n}:=\tau(H_n)=\frac{1}{\lambda_n}$, this implies that
$$
H_{n}+\sigma_{_{n+1}}\big(K_{X_{n}}+\Delta_{n}\big) \sim_{\r} (\varphi_n)^*\Big(H_{n+1}+\sigma_{_{n+1}}\big(K_{X_{n+1}}+\Delta_{n+1}\big)\Big) +a(\sigma_{_{n+1}}-\tau_{_n})E.
$$
Since $\tau_{_n}\leq \sigma_{_{n}}\leq \sigma_{_{n+1}}$, 
this shows that $H_{n}+\sigma(H_{n+1})\big(K_{X_{n}}+\Delta_{n}\big)$ is pseudo-effective, and thus $\sigma_{_{n+1}}= \sigma_{_n}$.
\end{proof}

\begin{thm} \label{thm:ample}
Let $(X,\Delta)$ be a $\q$-factorial klt pair, and suppose that $(K_X+\Delta)\notin \eff(X)$.
Let $H$ be an effective ample $\r$-divisor  on $X$ such that $K_X+\Delta+H$ is nef and klt.
Then the following holds.
\begin{enumerate}
\item Any $(K_X+\Delta)$-MMP with scaling of $H$ yields a birational map $\varphi:X\map X'$ and a Mori Fiber Space $f:X'\to Y'$
such that 
$$
\big(H+\sigma(H)(K_X+\Delta)\big)\cdot \varphi_{num}^*(\xi)=0,
$$
where $\xi\in N_1(X')$ is the class  of any curve on a fiber of $f$.
\item If $\Delta$ and $H$ are $\q$-divisors, then $\sigma(X,\Delta,H)\in \q$.
\end{enumerate}
\end{thm}

\begin{proof}
Let $\big\{(X_i, \varphi_i, \Delta_i, H_i, \lambda_i, R_i)\big\}$ be a $(K_X+\Delta)$-MMP with scaling of $H$
as described in \ref{MMP_with_scaling}.
By Theorem~\ref{mmp_w_scaling_terminates}, it terminates with  a Mori Fiber Space $f:X_n\to Y'$ associated to a  $(K_{X_n}+\Delta_n)$-negative 
extremal ray $R_n\subset \nec(X_n)$ such that $(K_{X_n}+\Delta_n + \lambda_n H_n)\cdot R_n=0$.  
Equivalently, 
\begin{equation} \label{Rn}
	\big(H_{n}+\tau(H_n)(K_{X_{n}}+\Delta_{n})\big)\cdot R_n=0.
\end{equation}
Thus there is an $\r$-Cartier $\r$-divisor $D_{Y'}$ on $Y'$ such that
$H_{n}+\tau(H_n)\big(K_{X_{n}}+\Delta_{n}\big)\sim_{\r} f^*D_{Y'}$. 
This shows that $[H_{n}]+\tau(H_n)\big[K_{X_{n}}+\Delta_{n}\big]$
lies on the boundary of $\eff(X_n)$. 
Since $H_n$ is big, $[H_n]\in \inter\big(\eff(X_n)\big)$, and thus $\sigma(H_n)=\tau(H_n)>0$.
Moreover, by Lemma~\ref{process}, 
\begin{equation} \label{tau=sigma}
	\sigma(H)=\sigma(H_{n})=\tau(H_n).
\end{equation}

Set $X'=X_n$, $\varphi=\varphi_{n-1}\circ\cdots\circ \varphi_0:X\map X'$.
If $\xi$ is the class of any curve lying on a fiber of $f:X'\to Y'$, then $\xi \in R_n$. 
Equations~\eqref{Rn} and \eqref{tau=sigma}, together with the properties of the numerical pullback of curves $\varphi_{num}^*:N_1(X')\to N_1(X)$, imply that 
$\big(H+\sigma(H)(K_X+\Delta)\big)\cdot \varphi_{num}^*(\xi)=0$.

To prove (2), assuming that both $\Delta$ and $H$ are $\q$-divisors, 
we apply Theorem~\ref{rationality} to the pair $(X_n,\Delta_n)$ and $\q$-divisor $H_n$, observing \eqref{tau=sigma}.
\end{proof}

\begin{proof}[{Proof of Theorem~\ref{generalization}}]
Let the notation and assumptions be as in Theorem~\ref{generalization}. We also assume that $K_X+\Delta\not \in \eff(X)$. Set
$$
W=\NE +\overline{\sum_{[C]\in \Sigma}\r_{\geq 0}[C]}.
$$
First we show that  ${W}\subset \NE +\nm$.
Let $[C]\in \Sigma$.
Then there is a birational map $\varphi: X\map X'$ and a Mori Fiber Space $f:X'\to Y'$
obtained by a $(K_X+\Delta)$-MMP, and a curve $C'$ on a fiber of $f$ such that $[C]=\varphi_{num}^*\big([C']\big)$.
Since $\rho(X'/Y')=1$, $[C']$ is numerically proportional to the class of any curve on a fiber of $f$.
So we may assume that $C'$ moves in a dominating family of curves on $X'$ and that it
does not meet the indeterminacy locus of $\varphi^{-1}$.
We may also assume that $C$ is the birational transform of $C'$, and thus it 
moves in a dominating family of curves on $X$.
Thus  $\r_{\geq 0}[C]\subset \nm$.

Now suppose that  $W\neq \NE +\nm$. Then,
since these are both closed convex cones in $N_1(X)$, and none of them contains a line, there
exists an $\r$-divisor $D$ on $X$  such that:
\begin{enumerate}
\item[$\bullet$]  $D\cdot z_0=0$ for some $z_0\in \nm\setminus\{0\}$,
\item[$\bullet$]  $D\cdot z \geq 0$ for every $z\in \NE +\nm$, and
\item[$\bullet$]  $D\cdot z>0$ for every $z\in {W}\setminus \{0\}$.
\end{enumerate}

We claim that there exists an ample $\r$-divisor on $X$ of the form $D-a(K_X+\Delta)$ for some  $a>0$.
Indeed, by  Remark~\ref{big}, it is enough to prove that
$\inter\big(\nef(X)\big)\cap \big(\r_{>0}[D]+\r_{>0}[-(K_X+\Delta)]\big)\neq
\emptyset$. Assume otherwise. 
Then there exists a hyperplane in $N^1(X)$ separating these two cones. 
This means that there exists an element
$l\in N_1(X)\setminus\{0\}$ such that 
\begin{align}
	&\alpha\cdot l\geq 0\ \text{ for every } \ \alpha\in \nef(X),  \label{6.1} \\
	&D\cdot l\leq 0,\ \text{ and } \notag \\
	-&(K_X+\Delta)\cdot l\leq 0. \label{6.2}
\end{align}
Inequality \eqref{6.1} implies that $l\in \overline{NE}_1(X)\setminus\{0\}$.
Together with inequality \eqref{6.2},  this implies  that  $l\in {W}\setminus \{0\}$.
But this is impossible since 
$D\cdot z>0$ for every $z\in {W}\setminus {0}$. This proves the claim.

So we let $H$ be an effective ample $\r$-divisor on $X$ of the form $H=D-a(K_X+\Delta)$ for some $a>0$. 
Then $D=H+a(K_X+\Delta)$ is pseudo-effective because $D\cdot z\geq 0$ for every $z\in \nm$, and
$\nm$ is the dual cone to $\eff(X)$. 
Since $D\cdot z_0=0$ for some $z_0\in \nm\setminus\{0\}$, $[D]$ lies on the boundary of $\eff(X)$.
Since $H$ is ample, $[H]\in \inter\big(\eff(X)\big)$.
Thus $\sigma(H)=a$, and $D=H+\sigma(H)(K_X+\Delta)$.

Notice that $K_X+\Delta +\frac{1}{\tau(H)}H$ is nef. 
Since $H$ is ample, by replacing it with a suitable $\r$-linearly equivalent $\r$-divisor if necessary,
we may assume that $K_X+\Delta +\frac{1}{\tau(H)}H$ is klt.
By Theorem~\ref{thm:ample}, there is a birational
map $\varphi:X\map X'$ and a Mori Fiber Space $f:X'\to Y'$,
obtained by a $(K_X+\Delta)$-MMP with scaling of $\frac{1}{\tau(H)}H$,  such that 
\begin{equation} \label{xii}
	D\cdot  \varphi_{num}^*(\xi)=\big(H+\sigma(H)(K_X+\Delta)\big)\cdot \varphi_{num}^*(\xi)=0
\end{equation}
for any class $\xi$ of a curve lying on a fiber of $f$. 
Let $C'\subset X'$ be a curve avoiding the indeterminacy locus of $\varphi^{-1}$, and let $C\subset X$
be the birational transform of $C'$.  Then $[C]=\varphi_{num}^*\big([C']\big)$.
Thus $[C]\in\Sigma\subset W$. But this is impossible since, by construction, 
$D\cdot z>0$ for every $z\in {W}\setminus {0}$.
This proves the theorem.
\end{proof}

\begin{rem}\label{exposed}
Let $N$ be a closed convex cone of $\r^n$ and let $R$ be an extremal ray of $N$.
We say that $R$ is  \emph{exposed} if there exists a linear functional $D\in N^\vee$ such that 
$D\cdot z=0$ for every $z\in R$ and $D\cdot z>0$ for every $z\in N\setminus R$.
It follows from \cite{straszewicz} that $N$ is the closed convex hull of its exposed extremal rays.
In particular, every extremal ray in $N$ is the limit of a sequence of exposed extremal rays.

Let $(X,\Delta)$ be a $\q$-factorial klt  pair, and $\Sigma\subset N_1(X)$ be as defined in Theorem~\ref{generalization}.
Let $R$ be a $(K_X+\Delta)$-negative exposed extremal ray of $ \NE +\nm$.
Then $R\subset \nm$, and the proof of Theorem~\ref{generalization} shows that $R=\r_{\geq 0}[C]$ for some $[C]\in \Sigma$.
\end{rem}

\begin{proof}[{Proof of Corollary~\ref{log_fano}}]
By \cite[Corollary 1.3.4]{bchm}, when $-(K_X+\Delta)$ is ample, $\nm$ is a rational polyhedral cone.
We refer also to \cite{barkowski_cone} for simple proof of this result
when $\Delta=0$ and $X$ is a smooth Fano threefold.
In particular, $\nm$ has finitely many extremal rays, all of which are exposed.
The corollary then follows from Remark~\ref{exposed}.
\end{proof}


\section{The cone of nef curves of terminal threefolds}
\label{end}

In this section we prove Theorem~\ref{main_thm}.

Recall that if $X$ is a projective $\q$-factorial terminal threefold, and $X'$ is 
obtained from $X$ by a MMP, then $X'$ is also $\q$-factorial and terminal.

\begin{lemma}\label{boundedness}
There exists a constant $N$ such that the following holds.
Let $X$ be a  projective $\q$-factorial terminal threefold, and 
$f:X\to Y$ a Mori Fiber Space. 
Let $B$ be a subset of codimension at least $2$ in $X$.
Then there exists a proper curve $C\subset X\setminus B$ 
lying on a fiber of $f$ and such that $-K_X\cdot C\leq N$.
\end{lemma}

\begin{proof}
Recall that $-K_{X}$ is $f$-ample. 
By \cite[Corollary 5.18]{kollar_mori}, $X$ has only isolated singularities.
In particular, if $\dim(Y)>0$, then the general fiber $F$ of $f$ is smooth and $-K_{F}$ is ample.
If $Y$ is a surface, then $F$ is a smooth rational curve and $-K_X\cdot F=2$.
If $Y$ is a curve, then $F$ is a smooth Del Pezzo surface, and $B\cap F$ is finite.
By \cite[Theorem 1.10]{kollar_mori},  there is a rational curve $C\subset F\setminus B$ such that $-K_X\cdot C=-K_F\cdot C\leq 3$.
So we are left to deal with the case when $X$ is a $\q$-factorial terminal Fano threefold
with $\rho(X)=1$.

By \cite{kawamata_boundedness}, $\q$-factorial terminal Fano threefolds
with Picard number $1$ form a bounded family.
In particular, there exist universal constants $M$ and $D$ such that
the following holds.
For any  $\q$-factorial terminal Fano threefold $X$ with $\rho(X)=1$,
$-MK_X$ is a very ample Cartier divisor and $(-K_X)^3\leq D$.
An explicit bound for $M$ is given in \cite{kollar_effective_v_ampleness}.
Set $N=M^2D>3$.
Let $B\subset X$ be a subset of codimension at least $2$.
By intersecting two general members of the linear system $|-MK_X|$,
we obtain a proper curve $C\subset X\setminus B$ such that
$-K_X\cdot C = -K_X\cdot M^2(-K_X)^2\leq M^2D =N$.
\end{proof}

\begin{proof}[{Proof of Theorem~\ref{main_thm}}]

Given an ample $\r$-divisor $A$ on $X$, 
let $\Sigma_{A}\subset \Sigma$ be the set consisting of classes $[C]\in \Sigma$ such that $(K_X+ A)\cdot C<0$.
By Theorem~\ref{generalization}, 
$$
\NEb + \nm = \Big(\NEb + \nm \Big)_{_{(K_X+A)\geq 0}} +\overline{\sum_{[C]\in \Sigma_A}\r_{\geq 0}[C]}.
$$
In order to prove the theorem, we
must show that the set of rays $\big\{\r_{\geq 0}[C]\ |\ [C]\in \Sigma_{A}\big\}$ is finite.

Let $N$ be as in Lemma~\ref{boundedness}.
Let $[C]\in \Sigma_{A}$.
Then there is a birational map $\varphi: X\map X'$ and a Mori Fiber Space $f:X'\to Y'$
obtained by a $(K_X+\Delta)$-MMP, and a curve $C'$ on a fiber of $f$ such that $[C]=\varphi_{num}^*\big([C']\big)$.
Since $\rho(X'/Y')=1$, $[C']$ is numerically proportional to the class of any curve on a fiber of $f$.
So, by Lemma~\ref{boundedness},  we may assume that $C'$ avoids the
indeterminacy locus of $\varphi^{-1}$ and   $0<-K_{X'}\cdot C'\leq N$. 
We may also assume that $C$ is the birational transform of $C'$.
Thus $[C]$ is an effective integral class  and 
$0<-K_{X}\cdot C\leq N$.

Now notice that there are only finitely many integral classes of curves $C\subset X$ such that
$0<-K_X\cdot C\leq N$ and $(K_X+ A)\cdot C<0$. Indeed, they all satisfy
$A\cdot C<N$.
So we conclude that $\big\{\r_{\geq 0}[C]\ |\ [C]\in \Sigma_{A}\big\}$
is finite.
\end{proof}

\bibliographystyle{amsalpha}
\bibliography{carolina.bib}

\end{document}